\def\version{0.28}

\def\journal{TAMS}
\def\titlep{Triangular C$^{*}$-bialgebra defined as
the direct sum of matrix algebras}
\documentclass[11pt]{article}
\usepackage{graphicx,ifthen}
\usepackage{amssymb}
\usepackage{amsmath}

\font\germ=eufm10 at12pt

\def\goth#1{\hbox{\germ#1}}

\newcommand{\qed}{\hbox{\rule[-2pt]{3pt}{6pt}}}
\newcommand{\qedh}{\hfill\qed \\}

\newcommand{\vv}{\vspace{.3in}}

\newcommand{\vep}{\varepsilon}

\newtheorem{Thm}{Theorem}[section]

\newtheorem{rem}[Thm]{Remark}

\newtheorem{defi}[Thm]{Definition}
\newtheorem{lem}[Thm]{Lemma}

\newtheorem{fact}[Thm]{Fact}

\newcommand{\ww}{\vv\noindent}

\newcommand{\kn}{\Large\bf
$K\hspace{-.4cm} N$
\Large\bf\vv }

\def\cal#1{\mathcal #1}

\def\pr{{\it Proof.}\quad}

\def\co#1{{\cal O}_{#1}}
\def\disp#1{{\displaystyle #1}}
\def\brl{branching law}
\def\bfsnl{{\rm BFS}_{N}(\Lambda)}

\addtocounter{footnote}{1}
\def\sftt#1{
\setcounter{equation}{0}
\addtocounter{footnote}{1}
\section{#1}
}

\def\ssft#1{\subsection{#1}}

\def\sssft#1{\subsubsection{#1}}

\begin{document}
%
% Personal data
%
\def\autherp{Katsunori Kawamura}
\def\emailp{e-mail: kawamura@kurims.kyoto-u.ac.jp.}
\def\addressp{{\small {\it College of Science and Engineering,
 Ritsumeikan University,}}\\
{\small {\it 1-1-1 Noji Higashi, Kusatsu, Shiga 525-8577, Japan}}
}

\def\infw{\Lambda^{\frac{\infty}{2}}V}
\def\zhalfs{{\bf Z}+\frac{1}{2}}
\def\ems{\emptyset}
\def\pmvac{|{\rm vac}\!\!>\!\! _{\pm}}
\def\vac{|{\rm vac}\rangle _{+}}
\def\dvac{|{\rm vac}\rangle _{-}}
\def\ovac{|0\rangle}
\def\tovac{|\tilde{0}\rangle}
\def\expt#1{\langle #1\rangle}
\def\zph{{\bf Z}_{+/2}}
\def\zmh{{\bf Z}_{-/2}}
\def\brl{branching law}
\def\bfsnl{{\rm BFS}_{N}(\Lambda)}
\def\scm#1{S({\bf C}^{N})^{\otimes #1}}
\def\mqb{\{(M_{i},q_{i},B_{i})\}_{i=1}^{N}}
\def\zhalf{\mbox{${\bf Z}+\frac{1}{2}$}}
\def\zmha{\mbox{${\bf Z}_{\leq 0}-\frac{1}{2}$}}
\newcommand{\mline}{\noindent
\thicklines
\setlength{\unitlength}{.1mm}
\begin{picture}(1000,5)
\put(0,0){\line(1,0){1250}}
\end{picture}
\par
 }
\def\ptimes{\otimes_{\varphi}}
\def\delp{\Delta_{\varphi}}
\def\delps{\Delta_{\varphi^{*}}}
\def\gamp{\Gamma_{\varphi}}
\def\gamps{\Gamma_{\varphi^{*}}}
\def\sem{{\sf M}}
\def\hdelp{\hat{\Delta}_{\varphi}}
\def\tilco#1{\tilde{\co{#1}}}
\def\ndm#1{{\bf M}_{#1}(\{0,1\})}
\def\cdm#1{{\cal M}_{#1}(\{0,1\})}
\def\tndm#1{\tilde{{\bf M}}_{#1}(\{0,1\})}
\def\sck{{\sf CK}_{*}}
\def\hdel{\hat{\Delta}}
% Boldfont
\def\ba{\mbox{\boldmath$a$}}
\def\bb{\mbox{\boldmath$b$}}
\def\bc{\mbox{\boldmath$c$}}
\def\bd{\mbox{\boldmath$d$}}
\def\be{\mbox{\boldmath$e$}}
\def\bk{\mbox{\boldmath$k$}}
\def\bm{\mbox{\boldmath$m$}}
\def\bp{\mbox{\boldmath$p$}}
\def\bq{\mbox{\boldmath$q$}}
\def\bu{\mbox{\boldmath$u$}}
\def\bv{\mbox{\boldmath$v$}}
\def\bw{\mbox{\boldmath$w$}}
\def\bx{\mbox{\boldmath$x$}}
\def\by{\mbox{\boldmath$y$}}
\def\bz{\mbox{\boldmath$z$}}
\def\bomega{\mbox{\boldmath$\omega$}}
\def\N{{\bf N}}
\def\lxm{L_{2}(X,\mu)}
\def\qtimes{\otimes_{\tilde{\varphi}}}
\def\ul#1{\underline{#1}}
%%%%%%%%%%%%%%%%%%%%%%%%%%%%%%%%%%%%%%%%%%%%%%%%
\def\titlepage{%\vspace{-4cm}

\noindent
{\bf 
\noindent
\thicklines
\setlength{\unitlength}{.1mm}
\begin{picture}(1000,0)(0,-300)
\put(0,0){\kn \knn\, for \journal\, Ver.\version}
\put(0,-50){\today}
\end{picture}
}
\vspace{-2.3cm}
\quad\\
{\small file: \textsf{tit01.txt,\, J1.tex}
\footnote{
${\displaystyle
\mbox{directory: \textsf{\fileplace}, 
file: \textsf{\incfile},\, from \startdate}}$}}
\quad\\
\framebox{
%\noindent
\begin{tabular}{ll}
\textsf{Title:} &
\begin{minipage}[t]{4in}
\titlep
\end{minipage}
\\
\textsf{Author:} &\autherp
\end{tabular}
}
%\mline
{\footnotesize	
\tableofcontents }
}
\def\ngt{{\bf N}_{\geq 2}}
\def\ngti{{\bf N}_{\geq 2}^{\infty}}

%%%%%%%%%%%%%%%%%%%%%%%%%%%%%%%%%%%%%%%%%%%%%%%%
%
%%%%%%%%% Cut from here %%%%%%%%%%
%\input comm.txt
%%%%%%%%% End of Cut %%%%%%%%%
%
%
\setcounter{section}{0}
\setcounter{footnote}{0}
\setcounter{page}{1}
\pagestyle{plain}

%
% Title
%
%%%%%%%%%%%%%%%%%%%%%%%%%%%
\title{\titlep}
\author{\autherp\thanks{\emailp}
\\
\addressp}
\date{}
\maketitle

%
% Abstract
%
\begin{abstract}
Let $M_{*}({\bf C})$ denote
the C$^{*}$-algebra defined as the direct sum of 
all matrix algebras $\{M_{n}({\bf C}):n\geq 1\}$.
It is known that 
$M_{*}({\bf C})$ has a non-cocommutative comultiplication
$\Delta_{\varphi}$.
We show that the C$^{*}$-bialgebra $(M_{*}({\bf C}),\Delta_{\varphi})$
has a universal $R$-matrix $R$ such that
the quasi-cocommutative C$^{*}$-bialgebra 
$(M_{*}({\bf C}),\Delta_{\varphi},R)$ is triangular.
\end{abstract}

\noindent
{\bf Mathematics Subject Classifications (2000).} 16W35, 81R50, 46K10. \\
\\
{\bf Key words.} universal $R$-matrix, triangular C$^{*}$-bialgebra.

%%%%%%%%%%%%%%%%%%%%%%%%%%%%%%%%%%%%%%%%%%%%%%%%%%%
%
% Section 1
%
\sftt{Introduction}
\label{section:first}
The purpose of this paper
is to construct a new triangular C$^{*}$-bialgebra such that
its universal $R$-matrix is defined by using 
a certain  set of arithmetic transformations. 
In this section, we show our motivation,
definitions and our main theorem.

%%%%%%%%%%%%%%%%%%%%%%%%%%%%%%%%%%%%%%%%%%%%%%%%%%%%%%%%%%%%%%%%%%%%%%%%
%
% subsection 1.1
%
\ssft{Motivation}
\label{subsection:firstone}
In this subsection, we roughly explain our motivation 
and the background of this study.
Explicit mathematical definitions will 
be shown after $\S$ \ref{subsection:firsttwo}.

For $n\geq 2$,
let $M_{n}({\bf C})$ denote the C$^{*}$-algebra 
of all $n\times n$ matrices
and we define $M_{1}({\bf C})={\bf C}$
for convenience.
Define the C$^{*}$-algebra $M_{*}({\bf C})$
as the direct sum of $\{M_{n}({\bf C}):n\geq 1\}$:
%
% Equation 1.1
%
\begin{equation}
\label{eqn:cuntztwo}
M_{*}({\bf C})=M_{1}({\bf C})\oplus M_{2}({\bf C})
\oplus M_{3}({\bf C})\oplus\cdots.
\end{equation}
In $\S$ 6.3 of \cite{TS02}, 
we constructed a non-cocommutative comultiplication $\Delta_{\varphi}$ 
of $M_{*}({\bf C})$ such that $(M_{*}({\bf C}),\delp)$
is a C$^{*}$-subbialgebra of a certain C$^{*}$-bialgebra.
As a C$^{*}$-algebra, $M_{*}({\bf C})$
is almost trivial and there is no new property,
but the bialgebra structure is new,
which is not a deformation of a known cocommutative bialgebra

On the other hand, 
in the theory of quantum groups,
a universal $R$-matrix for a quasi-cocommutative bialgebras 
is important for applications to mathematical physics 
and low-dimensional topology \cite{Drinfeld,GRS,Jimbo,Kassel}.
Especially, quasi-triangular (or braided)
bialgebras generate solutions of Yang-Baxter equation. 
As a stronger property,
a triangular bialgebra was introduced by Drinfel'd \cite{Drinfeld}. 
In this case,
the tensor category of all representations of the bialgebra
is symmetric (\cite{Kassel}, XIII.6 Exercises 1).
See also \cite{CFW,GRS,VanDaeleVanKeer}.

Our interest is to find a universal $R$-matrix of 
$(M_{*}({\bf C}),\delp)$ in (\ref{eqn:cuntztwo})
if there exists.
In this paper, we construct a universal $R$-matrix $R$ 
of $(M_{*}({\bf C}),\delp)$
defined as a double infinite sequence of 
permutation matrices arising from 
certain arithmetic transformations
of quotients and residues of positive integers.
Furthermore,
we show that 
the quasi-cocommutative C$^{*}$-bialgebra  $(M_{*}({\bf C}),\delp,R)$
is triangular.

%%%%%%%%%%%%%%%%%%%%%%%%%%%%%%%%%%%%%%%%%%%%%%%%%%%
% 
% subsection 1.1
%
\ssft{Definitions}
\label{subsection:firsttwo}
In this subsection,
we recall definitions of C$^{*}$-bialgebra
and universal $R$-matrix \cite{TS20}.
At first,
we prepare terminologies about C$^{*}$-bialgebra according to \cite{KV,MNW}.
%%%%%%%%%%%%%%%%%%%%%%%%%%%%%%%%%%%%%%%%%%%%%%%%%%%%%
%
% subsubsection 1.2.1
%
\sssft{C$^{*}$-bialgebra}
\label{subsubsection:firsttwoone}
For a C$^{*}$-algebra $A$,
let $A^{''}$ denote the enveloping von Neumann algebra of $A$.
The {\it multiplier algebra} ${\cal M}(A)$ of $A$ is defined by
%
% Equation 1.2
%
\begin{equation}
\label{eqn:ma}
{\cal M}(A)\equiv 
\{a\in A^{''}:aA\subset A,\,
Aa\subset A\}.
\end{equation}
Then ${\cal M}(A)$ is a unital C$^{*}$-subalgebra of $A^{''}$.
Especially,
$A={\cal M}(A)$ if and only if $A$ is unital.
The algebra ${\cal M}(A)$ is the completion of $A$
with respect to the strict topology.

For two C$^{*}$-algebras $A$ and $B$,
let ${\rm Hom}(A,B)$ and $A\otimes B$ 
denote the set of all $*$-homomorphisms from $A$ to $B$ and 
the minimal C$^{*}$-tensor product of $A$ and $B$, respectively.
A $*$-homomorphism from $A$ to $B$ is not always extended to the map 
from ${\cal M}(A)$ to ${\cal M}(B)$.
If $f\in {\rm Hom}(A,B)$ is surjective and both $A$ and $B$
are separable, then $f$ is extended to a surjective 
$*$-homomorphism of ${\cal M}(A)$ onto ${\cal M}(B)$.
We state that $f\in {\rm Hom}(A,{\cal M}(B))$ 
is {\it nondegenerate} if $f (A)B$ is dense in $B$.
If both $A$ and $B$ are unital and $f$ is unital,
then $f$ is nondegenerate. 
For $f\in {\rm Hom}(A,{\cal M}(B))$,
if $f$ is nondegenerate, then 
$f$ is called a {\it morphism} from $A$ to $B$ \cite{W3}.
If $f$ is a nondegenerate $*$-homomorphism from $A$ to $B$,
then we can regard $f$ as a morphism from $A$ to $B$
by using the canonical embedding of $B$ into ${\cal M}(B)$.
Each morphism $f$ from $A$ to $B$ can be extended uniquely 
to a homomorphism $\tilde{f}$ from ${\cal M}(A)$ to ${\cal M}(B)$
such that $\tilde{f}(m) f(b)a = f(mb)a$ 
for $m\in {\cal M}(B), b\in B$, and $a \in A$.
If $f$ is injective, then so is $\tilde{f}$.

A pair $(A,\Delta)$ is a {\it C$^{*}$-bialgebra}
if $A$ is a C$^{*}$-algebra with $\Delta\in {\rm Hom}(A,{\cal M}(A\otimes A))$ 
such that $\Delta$ is nondegenerate and the following holds:
%
% Equation 1.3
%
\begin{equation}
\label{eqn:bialgebratwo}
(\Delta\otimes id)\circ \Delta=(id\otimes\Delta)\circ \Delta.
\end{equation}
We call $\Delta$ the {\it comultiplication} of $A$.
A C$^{*}$-bialgebra $(A,\Delta)$ is {\it counital}
if there exists $\vep\in {\rm Hom}(A,{\bf C})$ such that
%
% Equation 1.4
%
\begin{equation}
\label{eqn:counit}
(\vep\otimes id)\circ \Delta= id = (id\otimes \vep)\circ
\Delta.
\end{equation}
We call $\vep$ the {\it counit} of $A$ and write $(A,\Delta,\vep)$ as 
the counital C$^{*}$-bialgebra $(A,\Delta)$ with the counit $\vep$.
Remark that we do not assume $\Delta(A)\subset A\otimes A$.
Furthermore,
$A$ has no unit for a C$^{*}$-bialgebra $(A,\Delta)$ in general.

%%%%%%%%%%%%%%%%%%%%%%%%%%%%%%%%%%%%%%%%%%%%%%%%%%%%%
%
% subsubsection 1.2.2
%
\sssft{Universal $R$-matrix}
\label{subsubsection:firsttwotwo}
We recall a unitary universal $R$-matrix and the quasi-cocommutativity
for a C$^{*}$-bialgebra \cite{TS20}.
%
% Definition 1.2
%
\begin{defi}
\label{defi:rmat}
Let $(A,\Delta)$ be a  C$^{*}$-bialgebra.
\begin{enumerate}
%(i)
\item
The map $\tilde{\tau}_{A,A}$ from ${\cal M}(A\otimes A)$ 
to ${\cal M}(A\otimes A)$
is the extended flip defined as
%
% Equation 1.5
%
\begin{equation}
\label{eqn:tilde}
\tilde{\tau}_{A,A}(X)(x\otimes y)\equiv \tau_{A,A}(X(y\otimes x))
\quad(X\in {\cal M}(A\otimes A),\,x, y\in A)
\end{equation}
where $\tau_{A,A}$ denotes the flip of $A\otimes A$.
%(ii)
\item
The map $\Delta^{op}$ from $A$ to ${\cal M}(A\otimes A)$
defined as
%
% Equation 1.6
%
\begin{equation}
\label{eqn:opposite}
\Delta^{op}(x)\equiv \tilde{\tau}_{A,A}(\Delta(x))\quad(x\in A)
\end{equation}
is called the opposite comultiplication of $\Delta$.
%
%(iii)
\item
A C$^{*}$-bialgebra $(A,\Delta)$ is cocommutative
if $\Delta=\Delta^{op}$. 
%(iv)
\item
An element $R$ in ${\cal M}(A\otimes A)$
is called a (unitary) universal $R$-matrix of $(A,\Delta)$
if $R$ is a unitary and 
%
% Equation 1.7
%
\begin{equation}
\label{eqn:univ}
R\Delta(x)R^{*}=\Delta^{op}(x)\quad
(x\in A).
\end{equation}
In this case,
we state that $(A,\Delta)$ is quasi-cocommutative
(or almost cocommutative \cite{CP}).
\end{enumerate}
\end{defi}

\noindent
We write a quasi-cocommutative C$^{*}$-bialgebra $(A,\Delta)$ with
a universal $R$-matrix $R$ as $(A,\Delta,R)$.
If $A$ is unital,
then ${\cal M}(A\otimes A)=A\otimes A$ and
$\tilde{\tau}_{A,A}=\tau_{A,A}$.
In addition, if $(A,\Delta)$ is quasi-cocommutative with
a universal $R$-matrix $R$, then $R\in A\otimes A$.

Next, we introduce quasi-triangular and triangular C$^{*}$-bialgebra
according to \cite{Drinfeld}.
%
% Definition 1.2
%
\begin{defi}
\label{defi:quasi}
Let $(A,\Delta, R)$ be 
a quasi-cocommutative C$^{*}$-bialgebra.	
\begin{enumerate}
%(i)
\item
%(\cite{Kassel}, p173)
$(A,\Delta, R)$
is quasi-triangular  (or braided \cite{Kassel}) if the following holds:
%
% Equation 1.8
%
\begin{equation}
\label{eqn:delone}
(\Delta\otimes id)(R)=R_{13}R_{23},
\quad
(id\otimes \Delta)(R)=R_{13}R_{12}
\end{equation}
where we use the leg numbering notation \cite{BS}.
%(ii)
\item
$(A,\Delta, R)$ is triangular if 
$(A,\Delta, R)$ is quasi-triangular and 
the following holds:
%
% Equation 1.9
%
\begin{equation}
\label{eqn:tri}
R\,\tilde{\tau}_{A,A}(R)=I
\end{equation}
where $\tilde{\tau}_{A,A}$ is as in (\ref{eqn:tilde})
and $I$ denotes the unit of ${\cal M}(A\otimes A)$.
\end{enumerate}
\end{defi}
Since both $\Delta\otimes id$ and $id\otimes \Delta$
are nondegenerate,  (\ref{eqn:delone}) makes sense.
The equation (\ref{eqn:tri}) is written as ``$R^{12}R^{21}=1$" 
in \cite{Drinfeld}.
In Appendix \ref{section:appone},
we will show basic facts about quasi-triangular
 C$^{*}$-bialgebras.

%%%%%%%%%%%%%%%%%%%%%%%%%%%%%%%%%%%%%%%%%%%%%%%%%%%%%%%%
%
% subsubsection 1.2.3
%
\sssft{Direct product and direct sum of C$^{*}$-algebras}
\label{subsubsection:firsttwothree}
For an infinite set $\{A_{i}:i\in \Omega\}$ of C$^{*}$-algebras,
there are separate notions of direct sum and product which do not
coincide with the algebraic ones \cite{Blackadar2006}.
We define two C$^{*}$-algebras $\prod_{i\in\Omega} A_{i}$
and $\bigoplus_{i\in\Omega} A_{i}$ as follows:
%
% Equation 1.10, 1.11
%
\begin{eqnarray}
\label{eqn:prod}
\disp{\prod_{i\in\Omega} A_{i}}\equiv &
\disp{\{(a_{i}):\|(a_{i})\|\equiv \sup_{i}\|a_{i}\|<\infty\},}\\
\nonumber 
\\
\label{eqn:plus}
\disp{\bigoplus_{i\in\Omega} A_{i}}\equiv &
\disp{\{(a_{i}):\|(a_{i})\|\to 0\mbox{ as }i\to\infty\}}
\end{eqnarray}
in the sense that for every $\vep>0$ there are only finitely many
$i$ for which $\|a_{i}\|>\vep$.
We call $\prod_{i\in\Omega} A_{i}$ and $\bigoplus_{i\in\Omega} A_{i}$ 
the {\it direct product} and the {\it direct sum} of $A_{i}$'s,
respectively. 
The algebra $\bigoplus_{i\in\Omega} A_{i}$ 
is a closed two-sided ideal of $\prod_{i\in\Omega} A_{i}$.
The algebraic direct sum $\oplus_{alg}\{A_{i}:i\in\Omega\}$ 
is a dense $*$-subalgebra of $\oplus\{A_{i}:i\in\Omega\}$.
Since ${\cal M}(\oplus_{i\in \Omega}A_{i})
\cong \prod_{i\in \Omega}{\cal M}(A_{i})$ (\cite{Blackadar2006}, II.8.1.3),
if $A_{i}$ is unital for each $i$, then
%
% Equation 1.12
%
\begin{equation}
\label{eqn:multiplus}
{\cal M}\Bigl(\bigoplus_{i\in \Omega}A_{i}\Bigr)\cong \prod_{i\in \Omega}A_{i}.
\end{equation}

%%%%%%%%%%%%%%%%%%%%%%%%%%%%%%%%%%%%%%%%%%%%%%%%%%%%%%%%%%
%
% subsection 1.3
%
\ssft{C$^{*}$-bialgebra $(M_{*}({\bf C}),\delp)$}
\label{subsection:firstthree}
In this subsection,
we recall the C$^{*}$-bialgebra $(M_{*}({\bf C}),\delp)$ \cite{TS02}.
Let $M_{*}({\bf C})$ be as in (\ref{eqn:cuntztwo}) and 
let $\{E^{(n)}_{i,j}\}$ denote the set of standard matrix 
units of $M_{n}({\bf C})$.
For $n,m\geq 1$,
define 
$\varphi_{n,m}\in 
{\rm Hom}(M_{nm}({\bf C}),M_{n}({\bf C})\otimes M_{m}({\bf C}))$ by
%
% Equation 1.13
%
\begin{equation}
\label{eqn:eij}
\varphi_{n,m}(E_{m(i-1)+j,m(i^{'}-1)+j^{'}}^{(nm)})=
E_{i,i^{'}}^{(n)}\otimes  E_{j,j^{'}}^{(m)}
\end{equation}
for $i,i^{'}\in\{1,\ldots,n\}$ and $j,j^{'}\in\{1,\ldots,m\}$.
By using $\{\varphi_{n,m}\}_{n,m\geq 1}$,
define two maps $\delp\in 
{\rm Hom}(M_{*}({\bf C}),\,M_{*}({\bf C})\otimes M_{*}({\bf C}))$
and $\vep\in {\rm Hom}(M_{*}({\bf C}),{\bf C})$ by
%
% Equation 1.14
%
\begin{equation}
\label{eqn:delp}
\delp(x)\equiv \sum_{m,l:\,ml=n}\varphi_{m,l}(x)\quad
\mbox{ when }x\in M_{n}({\bf C}),
\end{equation}
%
%
% Equation 1.15
%
\begin{equation}
\label{eqn:counittwob}
\vep(x)\equiv 0\quad\mbox{ when }x\in \oplus \{M_{n}({\bf C}):n\geq 2\},
\quad \vep(x)\equiv x\quad\mbox{ when }x\in M_{1}({\bf C}).
\end{equation}
Then $(M_{*}({\bf C}),\delp,\vep)$ is a 
counital C$^{*}$-bialgebra,
which is non-cocommutative. 
In fact,
%
% Equation 1.16
%
\begin{equation}
\label{eqn:noncocommutative}
\delp(E_{2,2}^{(6)})=
I_{1}\otimes E_{2,2}^{(6)}+E_{1,1}^{(2)}\otimes E_{2,2}^{(3)}
+E_{1,1}^{(3)}\otimes E_{2,2}^{(2)}+E_{2,2}^{(6)}\otimes I_{1}
\end{equation}
where $I_{1}$ denotes the unit of $M_{1}({\bf C})={\bf C}$.
The second and third terms show $\delp\ne \delp^{op}$.
Remark $\delp(M_{*}({\bf C}))\subset M_{*}({\bf C})\otimes M_{*}({\bf C})$.
The C$^{*}$-bialgebra 
$(M_{*}({\bf C}),\delp,\vep)$
satisfies the cancellation law
(\cite{TS02}, Proposition A.1), and 
it never has an antipode (\cite{TS02}, Lemma 3.2).

%%%%%%%%%%%%%%%%%%%%%%%%%%%%%%%%%%%%%%%%%%%%%%%%%%%%%%%%%%%%
%
% subsection 1.4
%
\ssft{Main theorem}
\label{subsection:firstfour}
In this subsection, we show our main theorem.
For $n\geq 1$, let $\{e^{(n)}_{i}\}_{i=1}^{n}$ denote 
the standard basis of the finite dimensional Hilbert space ${\bf C}^{n}$.
%
% Definition 1.3
%
\begin{defi}
\label{defi:rmatrixone}
Define the unitary transformation $R^{(n,m)}$ on 
${\bf C}^{n}\otimes {\bf C}^{m}$ by
%
% Equation 1.17
%
\begin{equation}
\label{eqn:urm}
R^{(n,m)}(e_{i}^{(n)}\otimes e_{j}^{(m)})\equiv 
e_{\ul{i}}^{(n)}\otimes e_{\ul{j}}^{(m)}
\end{equation}
for $(i,j)\in \{1,\ldots,n\}\times\{1,\ldots,m\}$
where the pair $(\ul{i},\ul{j})\in \{1,\ldots,n\}\times
\{1,\ldots,m\}$ is uniquely defined as the following integer equation:
%
% Equation 1.18
%
\begin{equation}
\label{eqn:mij}
m(i-1)+j=n(\ul{j}-1)+\ul{i}.
\end{equation}
\end{defi}
For example,
$m(i-1)+j$ divided by $n$ equals $\ul{j}-1$
with a remainder of $\ul{i}$ when $1\leq \ul{i}\leq n-1$.

By the natural identification 
${\rm End}_{{\bf C}}({\bf C}^{n}\otimes {\bf C}^{m})
\cong M_{n}({\bf C})\otimes M_{m}({\bf C})$,
$R^{(n,m)}$ is regarded as 
a unitary element in  $M_{n}({\bf C})\otimes M_{m}({\bf C})$ 
for each $n,m\geq 1$.
From (\ref{eqn:multiplus}),
${\cal M}(M_{*}({\bf C})\otimes M_{*}({\bf C}))
=\prod_{n,m\geq 1}M_{n}({\bf C})\otimes M_{m}({\bf C})$.
Hence the set $\{R^{(n,m)}\}_{n,m\geq 1}$ in (\ref{eqn:urm})
defines a unitary element $R$ in 
${\cal M}(M_{*}({\bf C})\otimes M_{*}({\bf C}))$:
%
% Equation 1.19
%
\begin{equation}
\label{eqn:urmtwo}
R\equiv (R^{(n,m)})_{n,m\geq 1}\in 
{\cal M}(M_{*}({\bf C})\otimes M_{*}({\bf C})).
\end{equation}
Then the main theorem is stated as follows.
%
% Theorem 1.4
%
\begin{Thm}
\label{Thm:main}
Let $(M_{*}({\bf C}),\delp)$ be as in $\S$ \ref{subsection:firstthree}.
\begin{enumerate}
%(i)
\item
The unitary $R$ in (\ref{eqn:urmtwo})
is a universal $R$-matrix of $(M_{*}({\bf C}),\delp)$.
%(ii)
\item
In addition to (i),
the quasi-cocommutative C$^{*}$-bialgebra $(M_{*}({\bf C}),\delp,R)$ 
is triangular.
\end{enumerate}
\end{Thm}

We discuss the meaning of $R$ in (\ref{eqn:urmtwo}) as follows.
%
% Remark 1.5
%
\begin{rem}
\label{rem:one}
{\rm
%\begin{enumerate}%(i) \item
From (\ref{eqn:mij}),
the operator $R^{(n,m)}$ in (\ref{eqn:urm}) is induced from
the arithmetic transformation $\chi_{n,m}$ defined as
%
% Equation 1.20
%
\begin{equation}
\label{eqn:chichib}
(i,j)\mapsto \chi_{n,m}(i,j)\equiv (\ul{i},\ul{j}).
\end{equation}
The map $\chi_{n,m}$ is a permutation
of the set $\{1,\ldots,n\}\times \{1,\ldots,m\}$.
For a given integer $N$ in $\{1,\ldots,nm\}$,
$(i,j),(\ul{i},\ul{j})\in\{1,\ldots,n\}\times \{1,\ldots,m\}$
are uniquely determined by
%
% Equation 1.21
%
\begin{equation}
\label{eqn:nem}
N=m(i-1)+j=n(\ul{j}-1)+\ul{i}.
\end{equation}
Hence 
both $(i,j)$ and $(\ul{i},\ul{j})$ are modifications of quotients 
and residues of $N$.
From this, $\chi_{n,m}$ means a transformation between
quotients and residues of a given integer
with respect to a pair of fixed integers $n$ and $m$.
For example,
%
% Equation 1.22
%
\begin{equation}
\label{eqn:values}
\chi_{2,3}(1,2)=(2,1),\quad
\chi_{2,3}(2,1)=(2,2).
\end{equation}
From this, $\chi_{2,3}\ne (\chi_{2,3})^{-1}$.
This implies $R^{2}\ne id$.
It is interesting that
the triangular structure of a bialgebra is induced from 
such arithmetic transformations.
%%(ii) \item \end{enumerate}
}
\end{rem}

In $\S$ \ref{section:second},
we will introduce locally triangular C$^{*}$-weakly coassociative system
as a generalization of $\{(M_{n}({\bf C}),\varphi_{n,m},R^{(n,m)}):n,m\geq 1\}$.
By using general statements in $\S$ \ref{section:second},
we will prove Theorem \ref{Thm:main} in $\S$ \ref{section:third}.

%%%%%%%%%%%%%%%%%%%%%%%%%%%%%%%%%%%%%%%%%%%%%%%%%%%%%%%%%%%%%%%%%%%%%%%%%%%%%%%%%%%
%
% Section 2
%
\sftt{C$^{*}$-weakly coassociative system}
\label{section:second}
In this section, we consider a general method
of construction of C$^{*}$-bialgebras
in order to prove Theorem \ref{Thm:main}.
%
% subsection 2.1
%
\ssft{Definitions}
\label{subsection:secondone}
According to $\S$ 3 in \cite{TS02},
we recall a general method to construct a C$^{*}$-bialgebra
from a set of C$^{*}$-algebras and $*$-homomorphisms among them.
We call $\sem$ a {\it monoid} if $\sem$ is a semigroup with unit.
%
% Definition 2.1
% 
\begin{defi}
\label{defi:axiom}
Let $\sem$ be a monoid with the unit $e$.
A data $\{(A_{a},\varphi_{a,b}):a,b\in \sem\}$
is a C$^{*}$-weakly coassociative system (= C$^{*}$-WCS) over $\sem$ if 
$A_{a}$ is a unital C$^{*}$-algebra for $a\in \sem$
and $\varphi_{a,b}$ is a unital $*$-homomorphism
from $A_{ab}$ to $A_{a}\otimes A_{b}$
for $a,b\in \sem$ such that
\begin{enumerate}
%(i)
\item
for all $a,b,c\in \sem$, the following holds:
%
% Equation 2.1
%
\begin{equation}
\label{eqn:wcs}
(id_{a}\otimes \varphi_{b,c})\circ \varphi_{a,bc}
=(\varphi_{a,b}\otimes id_{c})\circ \varphi_{ab,c}
\end{equation}
where $id_{x}$ denotes the identity map on $A_{x}$ for $x=a,c$,
%(ii)
\item
there exists a counit $\vep_{e}$ of $A_{e}$ 
such that $(A_{e},\varphi_{e,e},\vep_{e})$ 
is a counital C$^{*}$-bialgebra,
%(iii)
\item
$\varphi_{e,a}(x)=I_{e}\otimes x$ and
$\varphi_{a,e}(x)=x\otimes I_{e}$ for $x\in A_{a}$ and $a\in \sem$.
\end{enumerate}
\end{defi}

\noindent 
From this definition, the following holds.
%
% Theorem 2.2
% 
\begin{Thm}
\label{Thm:mainthree}
(\cite{TS02}, Theorem 3.1)	
Let $\{(A_{a},\varphi_{a,b}):a,b\in \sem\}$ be a C$^{*}$-WCS 
over a monoid $\sem$.
Assume that $\sem$ satisfies 
%
% Equation 2.2
%
\begin{equation}
\label{eqn:finiteness}
\#{\cal N}_{a}<\infty \mbox{ for each }a\in \sem
\end{equation}
where ${\cal N}_{a}\equiv\{(b,c)\in \sem\times \sem:\,bc=a\}$.
Define C$^{*}$-algebras 
%
% Equation 2.3
%
\begin{equation}
\label{eqn:astar}
A_{*}\equiv  \oplus \{A_{a}:a\in \sem\},\quad
C_{a}\equiv 
\oplus \{A_{b}\otimes A_{c}:(b,c)\in {\cal N}_{a}\}
\quad (a\in\sem),
\end{equation}
and define $*$-homomorphisms $\Delta^{(a)}_{\varphi}\in{\rm Hom}(A_{a},C_{a})$,
$\Delta_{\varphi}
\in {\rm Hom}(A_{*}, A_{*}\otimes A_{*})$ and
$\vep\in {\rm Hom}(A_{*},{\bf C})$ by 
%
% Equation 2.4
%
\begin{equation}
\label{eqn:cua}
\Delta^{(a)}_{\varphi}(x)\equiv \sum_{(b,c)\in {\cal N}_{a}}
\varphi_{b,c}(x)\quad(x\in A_{a}),\quad 
\Delta_{\varphi}\equiv \oplus\{\Delta_{\varphi}^{(a)}:a\in \sem\},
\end{equation}
%
% Equation 2.5
%
\begin{equation}
\label{eqn:counittwo}
\vep(x)\equiv 
\left\{
\begin{array}{cl}
0\quad &\mbox{ when }x\in \oplus \{A_{a}:a\in \sem\setminus \{e\}\},\\
\\
\vep_{e}(x)\quad&\mbox{ when }x\in A_{e}.
\end{array}
\right.
\end{equation}
Then $(A_{*},\delp,\vep)$ is a counital C$^{*}$-bialgebra.
\end{Thm}

\noindent
We call $(A_{*},\Delta_{\varphi},\vep)$ in 
Theorem \ref{Thm:mainthree} by a (counital)
{\it C$^{*}$-bialgebra} associated with 
$\{(A_{a},\varphi_{a,b}):a,b\in \sem\}$.
In this paper, we always assume the condition (\ref{eqn:finiteness}).

In this subsection, we do not assume that $\sem$ is abelian. 
For example, we constructed a C$^{*}$-WCS over a 
non-abelian monoid in \cite{TS05}.

%%%%%%%%%%%%%%%%%%%%%%%%%%%%%%%%%%%%%%%%%%%%%%%%%%%%%
%
% subsection 2.2
%
\ssft{Locally triangular C$^{*}$-weakly coassociative system}
\label{subsection:secondtwo}
In addition to $\S$ \ref{subsection:secondone},
we introduce locally triangular C$^{*}$-weakly coassociative system (=C$^{*}$-WCS) 
in this subsection.
%
% Definition 2.3
%
\begin{defi}
\label{defi:qqq}
Let 
$\{(A_{a},\varphi_{a,b}):a,b\in\sem\}$
be a  C$^{*}$-WCS.
\begin{enumerate}
%(i)
\item
For $a,b\in\sem$,
define $\varphi_{a,b}^{op}\in {\rm Hom}(A_{ab},A_{b}\otimes A_{a})$ by
%
% Equation 2.6
%
\begin{equation}
\label{eqn:vnm}
\varphi_{a,b}^{op}\equiv \tau_{a,b}\circ \varphi_{a,b}
\end{equation}
where $\tau_{a,b}$ denotes the flip from $A_{a}\otimes A_{b}$
to $A_{b}\otimes A_{a}$.
%(ii)
\item
$\{(A_{a},\varphi_{a,b}):a,b\in\sem\}$
is locally quasi-cocommutative
if there exists $\{R^{(a,b)}:a,b\in\sem\}$
such that $R^{(a,b)}$ is a unitary in
$A_{a}\otimes A_{b}$ 
and
%
% Equation 2.7
%
\begin{equation}
\label{eqn:rabv}
R^{(a,b)}\varphi_{a,b}(x)(R^{(a,b)})^{*}=\varphi_{b,a}^{op}(x)\quad
(x\in A_{ab})
\end{equation}
for each $a,b\in\sem$.
In this case,
we call
$\{(A_{a},\varphi_{a,b},R^{(a,b)}):a,b\in\sem\}$
a locally  quasi-cocommutative C$^{*}$-WCS. 
%(iii)
\item
A locally  quasi-cocommutative C$^{*}$-WCS 
$\{(A_{a},\varphi_{a,b},R^{(a,b)}):a,b\in\sem\}$
is locally  quasi-triangular
if 
the following holds:
%
% Equation 2.8, 2.9
%
\begin{eqnarray}
\label{eqn:varphiab}
(\varphi_{a,b}\otimes id_{c})(R^{(ab,c)})=&R^{(a,c)}_{13}R^{(b,c)}_{23},\\
\nonumber
\\
\label{eqn:varphiabc}
(id_{a}\otimes \varphi_{b,c})(R^{(a,bc)})=&R^{(a,c)}_{13}R^{(a,b)}_{12}
\end{eqnarray}
for each $a,b,c\in\sem$.
%(iv)
\item
A locally  quasi-cocommutative C$^{*}$-WCS 
$\{(A_{a},\varphi_{a,b},R^{(a,b)}):a,b\in\sem\}$
is locally  triangular
if 
$\{(A_{a},\varphi_{a,b},R^{(a,b)}):a,b\in\sem\}$
is locally  quasi-triangular and the following holds:
%
% Equation 2.10
%
\begin{equation}
\label{eqn:tautau}
R^{(a,b)}\tau_{b,a} (R^{(b,a)})=I_{a}\otimes I_{b}\quad(a,b\in\sem)
\end{equation}
where $I_{x}$ denotes the unit of $A_{x}$ for $x=a,b$.
\end{enumerate}
\end{defi}

For a C$^{*}$-WCS 
$\{(A_{a},\varphi_{a,b}):a,b\in\sem\}$,
the following holds from (\ref{eqn:multiplus}):
%
% Equation 2.11
%
\begin{equation}
\label{eqn:directtwo}
{\cal M}(A_{*}\otimes A_{*})\cong \prod_{a,b\in\sem}A_{a}\otimes A_{b}.
\end{equation}
Hence we identify an element in ${\cal M}(A_{*}\otimes A_{*})$
with that in $\prod_{a,b\in\sem}A_{a}\otimes A_{b}$.

By Definition \ref{defi:qqq},
the following holds.
%
% Lemma 2.4
%
\begin{lem}
\label{lem:quasi}
Assume that a monoid $\sem$ is abelian.
\begin{enumerate}
%(i)
\item
If a C$^{*}$-WCS 
$\{(A_{a},\varphi_{a,b}):a,b\in\sem\}$,
is locally  quasi-cocommutative with respect to $\{R^{(a,b)}:a,b\in\sem\}$ in
(\ref{eqn:rabv}),
then the unitary 
$R\in {\cal M}(A_{*}\otimes A_{*})$ defined by
%
% Equation 2.12
%
\begin{equation}
\label{eqn:sum}
R\equiv (R^{(a,b)})_{a,b\in\sem}
\end{equation}
is a universal $R$-matrix of $(A_{*},\delp)$.
%(ii)
\item
If a locally  quasi-cocommutative C$^{*}$-WCS 
$\{(A_{a},\varphi_{a,b},R^{(a,b)}):a,b\in\sem\}$
is locally  quasi-triangular,
then 
$(A_{*},\delp,R)$ is quasi-triangular for $R$ in (\ref{eqn:sum}).
%(iii)
\item
If a locally  quasi-cocommutative C$^{*}$-WCS 
$\{(A_{a},\varphi_{a,b},R^{(a,b)}):a,b\in\sem\}$
is locally  triangular,
then $(A_{*},\delp,R)$ is triangular for $R$ in (\ref{eqn:sum}).
\end{enumerate}
\end{lem}
%
% Proof
%
\pr
(i)
Let $a\in\sem$ and $x\in A_{a}$.
From (\ref{eqn:cua}),
%
% Equation 2.13
%
\begin{equation}
\label{eqn:relone}
R\delp(x)R^{*}
=
R\delp^{(a)}(x)R^{*}
=
\disp{\sum_{b,c;\,bc=a}R\varphi_{b,c}(x)R^{*}.}
\end{equation}
From (\ref{eqn:rabv}),
%
% Equation 2.14
%
\begin{equation}
\label{eqn:reltwo}
R\varphi_{b,c}(x)R^{*}=R^{(b,c)}\varphi_{b,c}(x)(R^{(b,c)})^{*}
=\varphi_{c,b}^{op}(x).
\end{equation}
From these and the assumption that $\sem$ is abelian,
%
% Equation 2.15
%
\begin{equation}
\label{eqn:relthree}
R\delp(x)R^{*}=
\sum_{b,c;\,bc=a}\varphi_{c,b}^{op}(x)
=\sum_{b,c;\,cb=a}\varphi_{c,b}^{op}(x).
\end{equation}
On the other hand,
%
% Equation 2.16
%
\begin{equation}
\label{eqn:commute}
\delp^{op}(x)=
\tilde{\tau}_{A_{*},A_{*}}(\delp^{(a)}(x))
=
\sum_{b,c;cb=a}\tau_{c,b}(\varphi_{c,b}(x))
=
\sum_{b,c;cb=a}\varphi_{c,b}^{op}(x).
\end{equation}
Hence 
$R\delp(x)R^{*}=\delp^{op}(x)$ for each $a\in\sem$ and $x\in A_{a}$.
Therefore the statement holds.

\noindent
(ii)
Let $a,b,c\in\sem$ and $z\in A_{a}\otimes A_{b}\otimes A_{c}$.
By (\ref{eqn:varphiab}),
%
% Equation 2.17
%
\begin{equation}
\label{eqn:seventeen}
(\delp\otimes id)(R)z
=
(\varphi_{a,b}\otimes id_{c})(R^{(ab,c)})z
=
R^{(a,c)}_{13}R^{(b,c)}_{23}z
=R_{13}R_{23}z.
\end{equation}
From this,
$(\delp\otimes id)(R)=R_{13}R_{23}$.
By the same token,
we can verify that 
$(id\otimes \delp)(R)=R_{13}R_{12}$.
Hence the statement holds.

\noindent
(iii)
Let $a,b\in\sem$ and $z\in A_{a}\otimes A_{b}$.
From (\ref{eqn:tautau}),
%
% Equation 2.18
%
\begin{equation}
\label{eqn:lastbb}
R\tilde{\tau}_{A_{*},A_{*}}(R)z=
R^{(a,b)}\tau_{b,a}(R^{(b,a)})z=z.
\end{equation}
This holds for each $a,b\in\sem$ and $z\in A_{a}\otimes A_{b}$.
Therefore $R\tilde{\tau}_{A_{*},A_{*}}(R)=I$.
Hence the statement holds.
\qedh

\noindent
We use the assumption that $\sem$ is abelian 
in the proof of Lemma \ref{lem:quasi}(i).

%%%%%%%%%%%%%%%%%%%%%%%%%%%%%%%%%%%%%%%%%%%%%
%
% Section 3
%
\sftt{Proof of Theorem \ref{Thm:main}}
\label{section:third}
We prove Theorem \ref{Thm:main} in this section.
We regard the set ${\bf N}\equiv \{1,2,3,\ldots\}$ 
of all positive integers
as a monoid with respect to
the multiplication. 
Then we see that $\{(M_{n}({\bf C}),\varphi_{n,m}):n,m\in {\bf N}\}$ 
in (\ref{eqn:eij})
is a C$^{*}$-WCS over the abelian monoid ${\bf N}$.
From Lemma \ref{lem:quasi},
it is sufficient for the proof of Theorem \ref{Thm:main} 
to show the following equations for $\{R^{(n,m)}:n,m\in {\bf N}\}$ in (\ref{eqn:urm}):
%
% Equation 3.1, 3.2, 3.3, 3.4
%
\begin{eqnarray}
R^{(n,m)}\varphi_{n,m}(x)(R^{(n,m)})^{*}
=&\varphi_{m,n}^{op}(x)\quad(x\in M_{nm}({\bf C})),
\label{eqn:braidzero}
\\
\nonumber
\\
\label{eqn:braidone}
(\varphi_{n,m}\otimes id_{l})(R^{(nm,l)})=&R_{13}^{(n,l)}R_{23}^{(m,l)},\\
\nonumber
\\
\label{eqn:braidtwo}
(id_{n}\otimes \varphi_{m,l})(R^{(n,ml)})=&R_{13}^{(n,l)}R_{12}^{(n,m)},\\
\nonumber
\\
\label{eqn:braidthree}
R^{(n,m)}\tau_{m,n}(R^{(m,n)})=&I_{n}\otimes I_{m}
\end{eqnarray}
for each $n,m,l\in {\bf N}$.
%%%%%%%%%%%%%%%%%%%%%%%%%%%%%%%%%%%%%%%%%%%
%
% subsection 3.1
%
\ssft{Proof of Theorem \ref{Thm:main}(i)}
\label{subsection:thirdone}
In this subsection,
we show (\ref{eqn:braidzero})
in order to prove Theorem \ref{Thm:main}(i).
We introduce several new symbols for convenience as follows:
Let  $F_{n}\equiv \{1,\ldots,n\}$ and
define the bijective map $\phi_{n,m}$ from $F_{n}\times F_{m}$
to $F_{nm}$ by 
%
% Equation 3.5
%
\begin{equation}
\label{eqn:nmphij}
\phi_{n,m}(i,j)\equiv m(i-1)+j\quad((i,j)
\in F_{n}\times F_{m}).
\end{equation}
Let $\{E^{(n)}_{i,j}:i,j\in F_{n}\}$ be as in $\S$ \ref{subsection:firstthree}.
For $i,j\in F_{n}$ and $k,l\in F_{m}$,
let 
%
% Equation 3.6
%
\begin{equation}
\label{eqn:fik}
{\bf E}_{(i,k),\,(j,l)}^{(n,m)}\equiv E_{i,j}^{(n)}\otimes E_{k,l}^{(m)}.
\end{equation}
%
%
% Lemma 3.1
%
\begin{lem}
\label{lem:addtion}
For $k,l\in F_{nm}$, the following holds:
%
% Equation 3.7, 3.8
%
\begin{eqnarray}
\label{eqn:hone}
R^{(n,m)}\varphi_{n,m}(E_{k,l}^{(nm)})(R^{(n,m)})^{*}
=&
{\bf E}_{\{\chi_{n,m}\circ \phi_{n,m}^{-1}\}(k),\,
\{\chi_{n,m}\circ\phi_{n,m}^{-1}\}(l)}^{(n,m)},\\
\nonumber 
\\
\label{eqn:eulbb}
\varphi_{m,n}^{op}(E_{k,l}^{(nm)})
=&{\bf E}^{(n,m)}_{\{\theta_{m,n}\circ \phi_{m,n}^{-1}\}(k),\,
\{\theta_{m,n}\circ \phi_{m,n}^{-1}\}(l)}
\end{eqnarray}
where $\chi_{n,m}$ is as in (\ref{eqn:chichib})
and
$\theta_{m,n}$ denotes the flip
from $F_{m}\times F_{n}$ to $F_{n}\times F_{m}$.
\end{lem}
%
% Proof
%
\pr
From (\ref{eqn:nmphij}) and (\ref{eqn:eij}),
%
% Equation 3.9
%
\begin{equation}
\label{eqn:crossd}
\varphi_{n,m}(E^{(nm)}_{k,l})
={\bf E}^{(n,m)}_{\phi_{n,m}^{-1}(k),\,\phi_{n,m}^{-1}(l)}
\quad(k,l\in F_{nm}).
\end{equation}
By definition,
$R^{(n,m)}$ is written as follows:
%
% Equation 3.10
%
\begin{equation}
\label{eqn:lastb}
R^{(n,m)}=\sum_{(i,j)\in F_{n}\times F_{m}} 
{\bf E}_{\chi_{n,m}(i,j),\,(i,j)}^{(n,m)}.
\end{equation}
From (\ref{eqn:crossd}) and (\ref{eqn:lastb}), (\ref{eqn:hone}) holds.

We see that
%
% Equation 3.11
%
\begin{equation}
\label{eqn:fikl}
\tau_{n,m}({\bf E}_{(i,k),\,(j,l)}^{(n,m)})=
E_{k,l}^{(m)}\otimes E_{i,j}^{(n)}
={\bf E}_{(k,i),\,(l,j)}^{(m,n)}.
\end{equation}
By (\ref{eqn:crossd}),
$\varphi_{m,n}(E^{(nm)}_{k,l})
={\bf E}^{(m,n)}_{\phi_{m,n}^{-1}(k),\,\phi_{m,n}^{-1}(l)}$.
From this and (\ref{eqn:fikl}), (\ref{eqn:eulbb}) holds.
\qedh

\noindent
The equation (\ref{eqn:lastb}) shows a representation of $R$ restricted
on the vector subspace $M_{n}({\bf C})\otimes M_{m}({\bf C})$
of $M_{*}({\bf C})\otimes M_{*}({\bf C})$.

\ww
{\it Proof of Theorem \ref{Thm:main}(i).}
We prove (\ref{eqn:braidzero}) for each $n,m\in {\bf N}$.
If it is done, then 
Theorem \ref{Thm:main}(i) holds 
from Lemma \ref{lem:quasi}(i).

For $\chi_{n,m}$ in (\ref{eqn:chichib}),
we see that
%
% Equation 3.12
%
\begin{equation}
\label{eqn:tauc}
\chi_{n,m}=\theta_{m,n}\circ \phi_{m,n}^{-1}\circ \phi_{n,m}
\end{equation}
where $\theta_{m,n}$ is as in Lemma \ref{lem:addtion}.
From this, we see that (\ref{eqn:hone}) equals (\ref{eqn:eulbb}).
Hence (\ref{eqn:braidzero}) is verified.
Therefore the statement holds.
\qedh

%%%%%%%%%%%%%%%%%%%%%%%%%%%%%%%%%%%%%%%%%%%%%%%%%
%
% subsection 3.2
%
\ssft{Proof of Theorem \ref{Thm:main}(ii)}
\label{subsection:thirdtwo}
In order to prove Theorem \ref{Thm:main}(ii),
we prove (\ref{eqn:braidone}-\ref{eqn:braidthree}).
Especially,
we show (\ref{eqn:braidone}) in
$\S$ \ref{subsubsection:thirdtwoone}
and $\S$ \ref{subsubsection:thirdtwotwo} step by step.
%%%%%%%%%%%%%%%%%%%%%%%%%%%%%%%%%%%%%%%%%%%%%%%%%%%%%
%
% subsubsection 3.2.1
%
\sssft{Proof of (\ref{eqn:braidone}) ---Step 1}
\label{subsubsection:thirdtwoone}
In this subsubsection,
we reduce (\ref{eqn:braidone}) to equations of maps on integers.
%
% Lemma 3.2
%
\begin{lem}
\label{lem:eqtwo}
Let $F_{n}$ and $\{E_{i,j}^{(n)}\}$ be as in $\S$ \ref{subsection:thirdone}.
For 
$i,a\in F_{n}$,
$j,b\in F_{m}$ and
$k,c\in F_{l}$,
let 
%
% Equation 3.13
%
\begin{equation}
\label{eqn:fijk}
{\bf E}_{(i,j,k),\,(a,b,c)}^{(n,m,l)}\equiv 
E_{i,a}^{(n)}\otimes 
E_{j,b}^{(m)}\otimes E_{k,c}^{(l)}.
\end{equation}
Then the following holds:
%
% Equation 3.14, 3.15
%
\begin{eqnarray}
\label{eqn:delpo}
(\varphi_{n,m}\otimes id_{l})(R^{(nm,l)})=&
\disp{\sum_{(a,b,c)\in F_{n}\times F_{m}\times F_{l}}
{\bf E}_{P(a,b,c),\,(a,b,c)}^{(n,m,l)},}\\
\nonumber\\
\label{eqn:rrr}
R_{13}^{(n,l)}R_{23}^{(m,l)}=&
\disp{\sum_{(a,b,c)\in F_{n}\times  F_{m}\times F_{l}}
{\bf E}_{Q(a,b,c),\,(a,b,c)}^{(n,m,l)}}
\end{eqnarray}
where $P$ and $Q$ are maps on $F_{n}\times F_{m}\times F_{l}$ defined by
%
% Equation 3.16, 3.17
%
\begin{eqnarray}
\label{eqn:pppa}
P\equiv   &
(\phi_{n,m}^{-1}\times id_l)
\circ \chi_{nm,l}\circ (\phi_{n,m}\times id_{l}),\\
\nonumber
\\
\label{eqn:ptwob}
Q\equiv  & (id_{n}\times \theta_{l,m})\circ (\chi_{n,l}\times id_{m})\circ
(id_{n}\times \theta_{m,l})\circ
(id_{n}\times \chi_{m,l})
\end{eqnarray}
where $id_{x}$ denotes the identity map on $F_{x}$ for $x=n,m,l$.
\end{lem}
%
% Proof
%
\pr
From (\ref{eqn:lastb}),\\
\\
$(\varphi_{n,m}\otimes id_{l})(R^{(nm,l)})$
\[\begin{array}{rl}
=&
\disp{\sum_{(t,k)\in F_{nm}\times F_{l}}
(\varphi_{n,m}\otimes id_{l})({\bf E}_{\chi_{nm,l}(t,k),\, (t,k)}^{(nm,l)})}\\
\\
=&
\disp{\sum_{(i,j,k)\in F_{n}\times F_{m}\times F_{l}}
(\varphi_{n,m}\otimes id_{l})
({\bf E}_{\chi_{nm,l}(\phi_{n,m}(i,j),k),\, (\phi_{n,m}(i,j),k)}^{(nm,l)})}.
\end{array}
\]
When $t= \phi_{n,m}(i,j)$,
\[
\begin{array}{rl}
(\varphi_{n,m}\otimes id_{l})({\bf E}_{\chi_{nm,l}(t,k), \,(t,k)}^{(nm,l)})
=&
\varphi_{n,m}(E^{(nm)}_{\ul{t},t})\otimes E^{(l)}_{\ul{k},k}\\
=&
{\bf E}^{(n,m)}_{\phi_{n,m}^{-1}(\ul{t}),\,\phi_{n,m}^{-1}(t)}\otimes E^{(l)}_{\ul{k},k}
\quad(\mbox{by (\ref{eqn:crossd})})\\
=&
{\bf E}^{(n,m)}_{\phi_{n,m}^{-1}(\ul{t}),\,(i,j)}\otimes E^{(l)}_{\ul{k},k}\\
=&
{\bf E}^{(n,m,l)}_{(\phi_{n,m}^{-1}(\ul{t}),\,\ul{k}),\,(i,j,k)}\\
\end{array}
\]
where $(\ul{t},\ul{k})=\chi_{nm,l}(t,k)$.
We see that
\[
\begin{array}{rl}
(\phi_{n,m}^{-1}(\ul{t}),\ul{k})
=&
\{(\phi_{n,m}^{-1}\times id_{l})\circ \chi_{nm,l}\}(t,k)\\
=&
\{(\phi_{n,m}^{-1}\times id_{l})\circ \chi_{nm,l}\circ (\phi_{n,m}\times id_{l})\}(i,j,k)
= P(i,j,k).
\end{array}
\]
Hence (\ref{eqn:delpo})  holds.

From (\ref{eqn:lastb}),
\[
\begin{array}{rl}
R_{13}^{(n,l)}R_{23}^{(m,l)}=&
\{id_{n}\otimes \tau_{l,m}\}(R^{(n,l)}\otimes id_{m})
\{id_{n}\otimes \tau_{m,l}\}
(id_{n}\otimes R^{(m,l)})\\
\\
=&
\disp{\sum_{(i,t)\in F_{n}\times F_{l}}\,
\sum_{(j,k)\in F_{m}\times F_{l}}
Y_{i,t,j,k}^{(n,m,l)}}
\end{array}
\]
where
%
% Equation 3.18
%
\begin{equation}
\label{eqn:yij}
Y_{i,t,j,k}^{(n,m,l)}\equiv 
\{id_{n}\otimes \tau_{l,m}\}
({\bf E}_{\chi_{n,l}(i,t),\,(i,t)}^{(n,l)}
\otimes id_{m}) 
\{id_{n}\otimes \tau_{m,l}\}
(id_{n}\otimes {\bf E}_{\chi_{m,l}(j,k),\,(j,k)}^{(m,l)}).
\end{equation}
Then we see that
%
% Equation 3.19
%
\begin{equation}
\label{eqn:nineteen}
Y_{i,t,j,k}^{(n,m,l)}
=
E_{\ul{i},i}^{(n)}\otimes E_{\ul{j},j}^{(m)}\otimes
E_{\ul{t},t}^{(l)}E_{\ul{k},k}^{(l)}
=\delta_{t,\ul{k}}\,{\bf E}_{(\ul{i},\ul{j},\ul{t}),\,(i,j,k)}^{(n,m,l)}
\end{equation}
where $(\ul{i},\ul{t})=\chi_{n,l}(i,t)$
and $(\ul{j},\ul{k})=\chi_{m,l}(j,k)$.
From this,
%
% Equation 3.20
%
\begin{equation}
\label{eqn:ijklmn}
R_{13}^{(n,l)}R_{23}^{(m,l)}=
\sum_{(i,j,k)\in F_{n}\times F_{m}\times F_{l}}
\left.{\bf E}_{(\ul{i},\ul{j},\ul{t}),\,(i,j,k)}^{(n,m,l)}\right|_{t=\ul{k}}
\end{equation}
When $t=\ul{k}$,
\[
\begin{array}{rl}
(\ul{i},\ul{j},\ul{t})
=&\{(id_{n}\times \theta_{l,m})\circ (\chi_{n,l}\times id_{m})\}(i,t,\ul{j})\\
=&\{(id_{n}\times \theta_{l,m})\circ (\chi_{n,l}\times id_{m})\}(i,\ul{k},\ul{j})\\
=&\{(id_{n}\times \theta_{l,m})\circ (\chi_{n,l}\times id_{m})\circ (id_{n}
\times \theta_{m,l})\circ (id_{n}\times \chi_{m,l})\}(i,j,k)\\
=& Q(i,j,k).
\end{array}
\]
Therefore (\ref{eqn:rrr}) holds.
\qedh

From Lemma \ref{lem:eqtwo},
it is sufficient for the proof of (\ref{eqn:braidone}) to show the equality 
$P=Q$ for two maps $P$ and $Q$ in (\ref{eqn:pppa}) and (\ref{eqn:ptwob}).

%%%%%%%%%%%%%%%%%%%%%%%%%%%%%%%%%%%%%%%%%%%%%%%%%%%%%%%
%
% subsubsection 3.2.2
%
\sssft{Proof of (\ref{eqn:braidone}) ---Step 2}
\label{subsubsection:thirdtwotwo}
In this subsubsection,
we prove equations of maps on integers in Lemma \ref{lem:eqtwo}.
%
% Lemma 3.3
%
\begin{lem}
\label{lem:equal}
For $P$ and $Q$ in (\ref{eqn:pppa}) and (\ref{eqn:ptwob}),
$P=Q$, that is,
the following diagram is commutative:

\def\diagram{
\put(400,500){$F_{n}\times F_{m}\times F_{l}$}
\put(100,380){$F_{nm}\times F_{l}$}
\put(100,230){$F_{nm}\times F_{l}$}
\put(400,100){$F_{n}\times F_{m}\times F_{l}$}
\put(800,400){$F_{n}\times F_{m}\times F_{l}$}
\put(800,300){$F_{n}\times F_{l}\times F_{m}$}
\put(800,200){$F_{n}\times F_{l}\times F_{m}$}
\thinlines
\put(350,470){\vector(-2,-1){80}}
\put(180,350){\vector(0,-1){70}}
\put(250,190){\vector(2,-1){80}}
\put(670,470){\vector(3,-1){100}}
\put(910,380){\vector(0,-1){30}}
\put(910,280){\vector(0,-1){30}}
\put(770,170){\vector(-3,-1){100}}
{\small
\put(120,470){$\phi_{n,m}\times id_{l}$}
\put(60,320){$\chi_{nm,l}$}
\put(120,120){$\phi_{n,m}^{-1}\times id_{l}$}
\put(730,470){$id_{n}\times \chi_{m,l}$}
\put(930,350){$id_{n}\times \theta_{m,l}$}
\put(930,250){$\chi_{n,l}\times id_{m}$}
\put(760,120){$id_{n}\times \theta_{l,m}$}
}
}

\noindent
\thicklines
%\framebox{
\setlength{\unitlength}{.1mm}
\begin{picture}(1200,480)(-50,70)
\put(0,0){\diagram}
\end{picture}
%}
\end{lem}
%
% Proof
%
\pr
Here we omit the symbol ``$\circ$" for simplicity of description.
For $\{\phi_{n,m}:n,m\in {\bf N}\}$ in (\ref{eqn:nmphij}),
the following holds:
%
% Equation 3.21
%
\begin{equation}
\label{eqn:sub}
\phi_{nm,l}(\phi_{n,m}\times id_{l})
=\phi_{n,ml}(id_{n}\times \phi_{m,l})\quad(n,m,l\in {\bf N}).
\end{equation}
From (\ref{eqn:sub}), we obtain
%
% Equation 3.22
%
\begin{equation}
\label{eqn:onea}
\phi_{nm,l}=
\phi_{n,ml}(id_{n}\times \phi_{m,l})(\phi_{n,m}\times id_{l})^{-1}.
\end{equation}
By the same token, we see that
%
% Equation 3.23, 3.24
%
\begin{eqnarray}
\label{eqn:twoa}
\phi_{n,ml}=&\phi_{n,lm}
=\phi_{nl,m}(\phi_{n,l}\times id_{m})(id_{n}\times \phi_{l,m})^{-1},\\
\nonumber
\\
\label{eqn:threea}
\phi_{nl,m}=&\phi_{ln,m}
=\phi_{l,nm}(id_{l}\times \phi_{n,m})(\phi_{l,n}\times id_{m})^{-1}.
\end{eqnarray}
Substituting  (\ref{eqn:threea}) into (\ref{eqn:twoa}),
and 
substituting it into (\ref{eqn:onea}),
%
% Equation 3.25
%
\begin{equation}
\label{eqn:nml}
\begin{array}{rl}
\phi_{nm,l}
=&\phi_{l,nm}(id_{l}\times \phi_{n,m})(\phi_{l,n}\times id_{m})^{-1}\\
&\times
 (\phi_{n,l}\times id_{m})(id_{n}\times \phi_{l,m})^{-1}
(id_{n}\times \phi_{m,l})(\phi_{n,m}\times id_{l})^{-1}\\
=&\phi_{l,nm}(id_{l}\times \phi_{n,m})(\theta_{n,l}\chi_{n,l}\times id_{m})
(id_{n}\times \theta_{m,l}\chi_{m,l})(\phi_{n,m}\times id_{l})^{-1}.
\end{array}
\end{equation}
Hence 
%
% Equation 3.26
%
\begin{equation}
\label{eqn:tauchi}
\theta_{nm,l}\chi_{nm,l}
=(id_{l}\times \phi_{n,m})(\theta_{n,l}\chi_{n,l}\times id_{m})
 (id_{n}\times \theta_{m,l}\chi_{m,l})(\phi_{n,m}\times id_{l})^{-1}.
\end{equation}
From this,
%
% Equation 3.27
%
\begin{equation}
\label{eqn:idl}
(id_{l}\times \phi_{n,m})^{-1}\theta_{nm,l}\chi_{nm,l}(\phi_{n,m}\times id_{l})
=(\theta_{n,l}\chi_{n,l}\times id_{m})(id_{n}\times \theta_{m,l}\chi_{m,l}).
\end{equation}
By multiplying $(id_{n}\times \theta_{l,m})(\theta_{l,n}\times id_{m})$
at both sides of (\ref{eqn:idl}) from the left,
%
% Equation 3.28
%
\begin{equation}
\label{eqn:idll}
\begin{array}{l}
(id_{n}\times \theta_{l,m})(\theta_{l,n}\times id_{m})
(id_{l}\times \phi_{n,m})^{-1}\theta_{nm,l}\chi_{nm,l}
(\phi_{n,m}\times id_{l})\qquad \qquad\\
\qquad\qquad =(id_{n}\times \theta_{l,m})(\chi_{n,l}\times id_{m})
(id_{n}\times \theta_{m,l}\chi_{m,l}).
\end{array}
\end{equation}
The R.H.S. of (\ref{eqn:idll}) is $Q$.
On the other hand, the L.H.S. of (\ref{eqn:idll}) is 
%
% Equation 3.29
%
\begin{equation}
\label{eqn:idt}
(id_{n}\times \theta_{l,m})(\theta_{l,n}\times id_{m})
\eta_{n,m,l}(\phi_{n,m}^{-1}\times id_{l})\chi_{nm,l}(\phi_{n,m}\times id_{l})
\end{equation}
where $\eta_{n,m,l}$ denotes the map from $F_{n}\times F_{m}\times F_{l}$
to $F_{l}\times F_{n}\times F_{m}$ defined as
$\eta_{n,m,l}(i,j,k)\equiv (k,i,j)$.
Since
$(id_{n}\times \theta_{l,m})(\theta_{l,n}\times id_{m})
\eta_{n,m,l}=id_{n}\times id_{m}\times id_{l}$,
the L.H.S. of (\ref{eqn:idll}) is $P$.
Hence the statement holds.
\qedh

\noindent
During initial phases of this study,
Lemma \ref{lem:equal} was forecasted by a computer experiment.
Essential pats of the proof of Lemma \ref{lem:equal}
are equations in (\ref{eqn:sub}).

%%%%%%%%%%%%%%%%%%%%%%%%%%%%%%%%%%%%%%%%%%%%%%%%%%%%%%
%
% subsubsection 3.2.3
%
\sssft{Proof of Theorem \ref{Thm:main}(ii)}
\label{subsubsection:thirdtwothree}
From Lemma \ref{lem:eqtwo} and Lemma \ref{lem:equal},
(\ref{eqn:braidone}) holds.
By the same token,
(\ref{eqn:braidtwo}) can be verified.
From Lemma \ref{lem:quasi}(ii),
the quasi-cocommutative C$^{*}$-bialgebra $(M_{*}({\bf C}),\delp,R)$
is quasi-triangular.

From (\ref{eqn:fikl}) and (\ref{eqn:lastb}),
\[\begin{array}{rl}
\tau_{n,m}(R^{(n,m)})
=&\disp{\sum_{(i,j)\in F_{n}\times F_{m}}
\tau_{n,m}({\bf E}_{\chi_{n,m}(i,j),\,(i,j)}^{(n,m)})}\\
\\
=&\disp{\sum_{(i,j)\in F_{n}\times F_{m}}
{\bf E}_{(\theta_{n,m}\chi_{n,m})(i,j),\,\theta_{n,m}(i,j)}^{(m,n)}}\\
\\
=&\disp{\sum_{(a,b)\in F_{n}\times F_{m}}
{\bf E}_{(\theta_{n,m}\chi_{n,m}\theta_{m,n})(b,a),\,(b,a)}^{(m,n)}.}\\
\end{array}
\]
From this and (\ref{eqn:lastb}),
\[
\begin{array}{rl}
R^{(n,m)}\tau_{m,n}(R^{(m,n)})
=&
\disp{
\sum_{(i,j),(b,a)\in F_{n}\times F_{m}}
{\bf E}_{\chi_{n,m}(i,j),\,(i,j)}^{(n,m)}
{\bf E}_{(\theta_{m,n}\chi_{m,n}\theta_{n,m})(b,a),\,(b,a)}^{(n,m)}}\\
\\
=&
\disp{
\sum_{(b,a)\in F_{n}\times F_{m}}
{\bf E}_{(\chi_{n,m}\theta_{m,n}\chi_{m,n}\theta_{n,m})(b,a),\,(b,a)}^{(n,m)}.}\\
\end{array}
\]
On the other hand,
$\chi_{n,m}\theta_{m,n}\chi_{m,n}\theta_{n,m}=id_{n}\times id_{m}$
from (\ref{eqn:tauc}).
Hence 
%
% Equation 3.30
%
\begin{equation}
\label{eqn:trial}
R^{(n,m)}\tau_{m,n}(R^{(m,n)})
=\sum_{(b,a)\in F_{n}\times F_{m}}{\bf E}_{(b,a),\,(b,a)}^{(n,m)}
=I_{n}\otimes I_{m}.
\end{equation}
Hence, (\ref{eqn:braidthree}) holds.
From this and Lemma \ref{lem:quasi}(iii),
the quasi-triangular C$^{*}$-bialgebra $(M_{*}({\bf C}),\delp,R)$
is triangular.
\qedh

\appendix

\section*{Appendix}

%%%%%%%%%%%%%%%%%%%%%%%%%%%%%%%%%%%%%%%%%%%%%%%%%%%
%
% Appendix A
%
\sftt{Basic facts about quasi-triangular C$^{*}$-bialgebras}
\label{section:appone}
In this section, we show basic facts 
about quasi-triangular C$^{*}$-bialgebras.
%
% Fact A.1
%
\begin{fact}
\label{fact:ybe}
Let $(A,\Delta,R)$ be a quasi-triangular C$^{*}$-bialgebra.
Then the following holds:
\begin{enumerate}
%(i)
\item
$R$ satisfies the Yang-Baxter equation
%
% Equation A.1
%
\begin{equation}
\label{eqn:ybe}
R_{12}R_{13}R_{23}=R_{23}R_{13}R_{12}.
\end{equation}
%
%(ii)
\item
If $(A,\Delta)$ has a counit $\vep$,
then
%
% Equation A.2
%
\begin{equation}
\label{eqn:vti}
(\vep\otimes id)(R)=id=(id\otimes \vep)(R)
\end{equation}
where $id$ denotes the unit of ${\cal M}(A)$.
%
%(iii)
\item
Let $({\cal H},\pi)$ be a nondegenerate representation
of the C$^{*}$-algebra $A$. Let $\Pi$
denote the extension  of $\pi\otimes \pi$ on ${\cal M}(A\otimes A)$
and let $T$ denote the flip on ${\cal H}\otimes {\cal H}$.
Define the unitary operator $C$ on ${\cal H}\otimes {\cal H}$ by
%
% Equation A.3
%
\begin{equation}
\label{eqn:athree}
C\equiv T\Pi(R).
\end{equation}
For $n\geq 3$,
let ${\cal H}^{\otimes n}$ denote
the $n$-times tensor power of ${\cal H}$.
For $1\leq i\leq n-1$,
let $C_{i}\equiv I_{{\cal H}}^{\otimes (i-1)}\otimes C\otimes I_{{\cal H}}^{\otimes (n-i)}$
where $I_{{\cal H}}$ denotes the identity map on ${\cal H}$.
Then
%
% Equation A.4
%
\begin{equation}
\label{eqn:bbb}
C_{i}C_{i+1}C_{i}=C_{i+1}C_{i}C_{i+1}.
\end{equation}
In addition,
if $(A,\Delta,R)$ is triangular,
then $C^{2}=I$.
\end{enumerate}
\end{fact}
%
% Proof 
%
\pr
Proofs of (i) and (ii) are given along with the proof 
of Theorem VIII.2.4 of \cite{Kassel}
which is modified to a C$^{*}$-bialgebra as follows:\\
(i)
\[
\begin{array}{rl}
R_{12}R_{13}R_{23}
=&R_{12}(\Delta\otimes id)(R)\quad(\mbox{by (\ref{eqn:delone})})\\
=&(\Delta^{op}\otimes id)(R)R_{12}\quad (\mbox{by (\ref{eqn:univ})})\\
=&(\tilde{\tau}_{A,A}\otimes id) ((\Delta\otimes id)(R))R_{12}\\
=&(\tilde{\tau}_{A,A}\otimes id) (R_{13}R_{23})R_{12}
\quad(\mbox{by (\ref{eqn:delone})})\\
=&
(\tilde{\tau}_{A,A}\otimes id) (R_{13})\cdot 
(\tilde{\tau}_{A,A}\otimes id) (R_{23})R_{12}\\
=&R_{23}R_{13}R_{12}.
\end{array}
\]

\noindent
(ii)
Since $\vep$ is nondegenerate,
it can be extended to the $*$-homomorphism $\tilde{\vep}$ from
${\cal M}(A)$ to ${\bf C}$ such that $\tilde{\vep}(I)=1$.
We write $\tilde{\vep}$ as $\vep$ here.
Since $(\vep\otimes id)\circ \Delta=id$,
\[
\begin{array}{rl}
R=&
\{(\vep\otimes id\otimes id)\circ (\Delta\otimes id)\}(R)\\
=&
(\vep\otimes id\otimes id)(R_{13}R_{23})
\qquad(\mbox{by (\ref{eqn:delone})})\\
=&
(\vep\otimes id\otimes id)(R_{13})\cdot
(\vep\otimes id\otimes id)(R_{23})\\
=&
(\vep\otimes id)(R)\cdot \vep(I)R.\\
\end{array}
\]
From this,
we obtain $(\vep\otimes id)(R)=id$
because $\vep(I)=1$ and $R$ is invertible.
By the same token,
we obtain $(id\otimes \vep)(R)=id$.

\noindent
(iii)
Assume $n=3$ and $i=1$.
Let $U\equiv \Pi(R)$.
Then
\[
\begin{array}{rl}
C_{1}C_{2}C_{1}
=&
T_{12}U_{12}
T_{23}U_{23}
T_{12}U_{12}\\
=&
T_{12}T_{23}T_{12}U_{23}U_{13}U_{12}\\
=&
T_{23}T_{12}T_{23}U_{12}U_{13}U_{23}\quad(\mbox{by (\ref{eqn:ybe})})\\
=&
T_{23}U_{23}T_{12}U_{12}T_{23}U_{23}\\
=& C_{2}C_{1}C_{2}\\
\end{array}
\]
where we use the leg numbering notations 
$T_{ij}$ and $U_{ij}$ on ${\cal H}^{\otimes 3}$
and $T_{12}T_{23}T_{12}=T_{23}T_{12}T_{23}$.
This implies (\ref{eqn:bbb}).

Assume that $(A,\Delta,R)$ is triangular.
For $a,b\in A$,
we see that 
$T\{(\pi\otimes \pi)(a\otimes b)\}T=(\pi\otimes \pi)(b\otimes a)$.
From this, 
%
% Equation A.5
%
\begin{equation}
\label{eqn:relfive}
T\Pi(R)T=\Pi(\tilde{\tau}_{A,A}(R)).
\end{equation}
From (\ref{eqn:relfive}) and (\ref{eqn:tri}),
%
% Equation A.6
%
\begin{equation}
\label{eqn:lastapp}
C^{2}= T\Pi(R)T\Pi(R)= \Pi(\tilde{\tau}_{A,A}(R)R)
= \Pi(I_{{\cal M}(A\otimes A)})=I.
\end{equation}
\qedh

In addition to Fact \ref{fact:ybe}(iii),
it is clear that 
$\{C_{i}\}_{i=1}^{n-1}$ satisfies $C_{i}C_{j}=C_{i}C_{j}$
for $i,j=1,\ldots,n-1$ when $|i-j|\geq 2$.
Therefore 
a nondegenerate representation of 
a quasi-triangular ({\it resp}. triangular) C$^{*}$-bialgebra
gives a unitary representation of 
the braid group $B_{n}$ (\cite{Kassel}, Lemma X.6.4)
({\it resp}. the symmetric group ${\goth S}_{n}$
 (\cite{Kassel}, $\S$ X.6.3)).

%\ww
%{\bf Acknowledgment:}
%The author would like to express his sincere thanks to Izumi Ojima 
%for his interest in this topic and raising the above question.

%%%%%%%%%%%%%%%%%%%%%%%%%%%%%%%%%%%%%%%%%%%%%%%%%%%%%%%%%%%%%%
%
% Reference 
%

%
\end{document}